\documentclass[12pt,a4paper]{article}

\usepackage[applemac]{inputenc}
\usepackage{amsfonts}
\usepackage{amsmath}
\usepackage{enumerate}

\long\def\tot#1\leb{}

\usepackage[colorlinks]{hyperref}
\newcommand{\R}{\mathbb R}
\newcommand{\N}{\mathbb N}
\newcommand{\Q}{\mathbb Q}
\newtheorem{definition}{Definition}
\newtheorem{theorem}{Theorem}
\newtheorem{lemma}{Lemma}
\newtheorem{corollary}{Corollary}

\def\Ord{{\cal O}}
\def\ord{{o}}

\newenvironment{proof}{\begin{trivlist}\item[]{\emph{Proof.}}}
               {\hfill$\Box$\end{trivlist}}

\newcommand*{\QEDB}{\hfill\ensuremath{\square}}%

\def\bcolor{\textcolor[rgb]{0.00,0.00,1.00}}

\title{
Asymptotically uniform functions:\\
a single hypothesis which solves\\
two old problems}
\author{Jean--Pierre Gabriel$^*$ and Jean--Paul Berrut\footnote{%
D{\'e}partement de Math{\'e}matiques,
Universit{\'e} de Fribourg, P{\'e}rolles, CH-1700 Fribourg, Switzerland;
{email: \tt jean-pierre.gabriel@unifr.ch, jean-paul.berrut@unifr.ch}}}

\begin{document}
\baselineskip15pt

\maketitle

\baselineskip12pt
{\abstract{The asymptotic study of a time-dependent function $f$ as the solution
of a differential equation often leads to the question of whether its derivative $\dot f$
vanishes at infinity. We show that a necessary and sufficient condition for this
is that $\dot f$ is what may be called {\sl asymptotically uniform.}
We generalize the result to higher order derivatives.
We further show that the same property for $f$ itself is also ne\-ces\-sa\-ry and sufficient
for its one-si\-ded improper integrals to exist. On the way, the article provides a
broad study of such asymptotically uniform functions.}}

\bigskip
{\small

\noindent {\em Math Subject Classification: }26A09, 26A12, 34A45

\smallskip

\noindent {\em Keywords: }asymptotically uniform function,
vanishing of a derivative at infinity,
vanishing of an integrand at infinity,
Hadamard's lemma,
Barb\u alat's lemma}

\baselineskip15pt
\section{Introduction}
When does the derivative of a function $f:\R_+\to \R$ tend to zero at infinity?
This certainly is the case when both $f$
and its derivative $\dot f := df/dt$ converge (see Section 2).
However,
it has been known for a long time that
this set of hy\-po\-the\-ses is too strong.


A similar situation occurs with (Riemann) integrals:
the fact that a function is integrable does not imply that
it tends to zero at infinity.

In the present work we shall 
address the following two 
questions.

\vskip1em
\noindent
$\bf{Differential\ case}$: Let $f:\R_+\rightarrow \R$ be differentiable with
$\lim_{t\rightarrow +\infty}f(t)=\alpha \in\R$.
Is there a necessary and sufficient condition for
$\dot f$ to converge to 0 at infinity?

\vskip0.5em
\noindent
$\bf{ Integral\ case}$: Let $g:\R_+\rightarrow \R$ be integrable
over every interval
$[0,t]$, with
$\lim_{t\rightarrow +\infty}$ \bcolor{ }$\int_0^t g(s)ds=\alpha \in\R$.
Does there exist a necessary and sufficient condition for
$g$ to converge to 0 at infinity?

\vskip1em
Since the domain of definition of all our functions will be
$\R_+$, we shall write 
$\infty$ instead of 
$+ \infty$.

\section{General discussion}
Historically, the dif\-fe\-ren\-tial case seems to have received more attention.
A simple result states that, if $f$ and $\dot f$ converge, then the limiting value of the
latter has to be 0.
This is an immediate consequence of the following lemma.

\begin{lemma}
Let $f:\R_+\rightarrow \R$ be differentiable with
$\lim_{t\rightarrow\infty}f(t)=\alpha \in \R$.
Then, for any sequence of intervals $([a_n,b_n])_{n\geq 1}$ with  $a_n<b_n$,
$a_n\uparrow\infty$ and $\inf_{n\geq 1}(b_n-a_n)>0$, there exists
$\xi_n\in (a_n,b_n)$ with $\lim_{n\rightarrow\infty}\dot f(\xi_n)=0$.
\end{lemma}
\begin{proof}
For every $n\geq 1$, the mean value theorem provides $a_n<\xi_n<b_n$ with
$f(b_n)-f(a_n)=\dot f(\xi_n)(b_n-a_n)$. Since $\inf_{n\geq 1} (b_n-a_n)$ is positive,
the convergence to $0$ of $f(b_n)-f(a_n)$ entails
$\lim_{n\rightarrow\infty}\dot f(\xi_n)=0$.
\end{proof}

Numerous converging functions have been given whose first de\-ri\-va\-tives
do not converge \cite[p.\ 40, and p.\ 45 for counterexamples in the integral case]{Gel-Olm}.

To our knowledge, 
apart from Thieme, who provided a necessary
and sufficient condition for the differential case \cite[Corollary A.17]{Thi},
only sufficient conditions have been given so far.
The uniform continuity of $\dot f$, res\-pec\-ti\-ve\-ly of $g$ for the
integral case, is one of the
best known and most frequently used among those sufficient conditions;
in control theory
\cite{Mul}, the 
result for $g$ is known as Barb\u alat's lemma \cite{Bar, Far-Weg}.


If one assumes continuity of $g$, the integral case is a consequence
of the differential case.
However,
in view of the existence of differentiable functions whose derivatives are not
integrable,
i.e., for which the fun\-da\-men\-tal theorem of calculus does not hold,
one has to treat the two ques\-tions separately:
recall for instance the existence of differentiable functions whose derivatives
are boun\-ded
but not integrable \cite[p.\ 88]{vRS}.


Already in the nineteenth century,
a considerable number of mathematicians
looked for sufficient conditions for the differential case.
As mentioned in \cite{Had2},
the first conditions for the va\-nishing of the derivative
seems to be studies by A.\ Kneser (1896)
about the asymptotic be\-ha\-viour of solutions $y$ of ordinary differential equations
$\ddot y = F(t,y)$, 
when $F$ possesses the same sign as $y$
\cite[p.\ 182]{Kne}. In his considerations on the asymptotics of a function
of the solution of a
system of ordinary differential equations, Hadamard \cite[p.\ 334]{Had1} appears to be
the first to have stated the result
independently of such equations; his supplementary assumption is the boundedness
of the second derivative.

\begin{theorem}[Hadamard's lemma \cite{Had1}]\label{th:Hadamard}
If the first $n$ derivatives $f^{(1)},$ $\ldots,f^{(n)}$ of  a function $f:\R_+\rightarrow \R$
are bounded, then the convergence of $f$ to a real number when $t\rightarrow\infty$
implies the convergence to $0$ of $f^{(1)},\ldots,f^{(n-1)}$.
\end{theorem}

\rm
Obviously unaware of \cite{Had1}, Littlewood \cite{Lit} gave an easier proof
of the result
for an arbitrarily differentiable function, which he deems ``fundamental''
for his study of power series.
In 1913, under somewhat stronger 
(but, as they stress, not necessary) conditions on $\ddot f$,
Hardy and Littlewood \cite{Har-Lit} initiated the study of bounds
on the first derivative
at infinity 
involving $f$ and $\ddot f$.
We give one of their results.
\begin{theorem}
Let $g$ and $h$ be positive non decreasing functions and $\ddot f$ be continuous. If
$f = \ord(g)$, or if $g = 1$ and $f\to \alpha$, and if $\ddot f = \Ord(h)$, then
$\dot f = \ord(\sqrt{gh})$.
\end{theorem}

In the same article,
the authors extended the result to the case of several derivatives; this
led them to the following version of Ha\-da\-mard's lemma (p.\ 423),
which merely requires the boundedness of
the last derivative.

\begin{theorem}
If $f$ has derivatives of all orders and if $f \to s$ at infinity,
then, if any derivative is bounded,
all preceding derivatives tend to zero, or, to state the matter in another way, if any
derivative does not tend to zero, no subsequent derivative is bounded.
\end{theorem}


\rm
A different proof of this result has been given by Coppel
\cite[p.\ 141]{Cop} in 1965.

In 1914, Landau \cite{Lan} originated a quantitative study of the problem by
obtaining bounds on $|f|$ and $|\ddot  f|$
which guarantee that $|\dot f(t)|$ lies below a certain value for sufficiently large $t$.
This result led
to an enormous body of research, exhibited in the book \cite{Mit-Pec-Fin},
which mentions it on its very first page.

Most results involve assumptions about $\ddot f$
for the differential case, on $\dot g$ for the integral one.
We rather look for
conditions
on $\dot f$ only 
that solve the differential case, resp.\ on $g$ only for the
integral one.
As mentioned already,
a class of functions introduced by Thieme
provides the solution to the differential case.
We show that it solves the integral case as well.
Furthermore, we perform a detailed study of this class, and
also give a more suggestive characterization of the latter.

The next section introduces the class of asymptotically uniform functions,
the subject of this manuscript, and shows that it solves the two problems
mentioned in the Introduction. Section 4 lists and discusses properties of
these functions. Section 5 studies some sufficient and/or necessary conditions for
a function itself to vanish at infinity. In Section 6 we show that
the class of asymptotically uniform functions may also be characterized as
those which can be written as the sum of a uniformly continuous
function and one which vanishes at infinity.
The last section extends
the differential case to functions with more derivatives



\section{Asymptotically uniform functions}
To our knowledge, for the integral case,
the usual sufficient conditions presuppose the continuity
of the integrand.
However, there exist asymptotically discontinuous functions, i.e., functions
discontinuous on every unbounded interval, 
which simultaneously have a convergent integral and vanish at infinity,
for instance on $\R^+$,
$$g(t) := \begin{cases}
0,&t\in\N,\\
e^{-t},&\hbox{elsewhere.}\\
\end{cases}$$
On every interval $[0,T]$, $T>0$, $g$
is Riemann integrable, as it has only finitely many discontinuities.
As $T\to\infty$, the integral converges to 1 and $g$ to zero.
For a necessary and sufficient condition,
we 
therefore include asymptotically discontinuous functions
in our considerations and introduce the following class, which
is equivalent to that proposed by Thieme. Additional
examples are provided in Section 4.

\begin{definition}
A function $f:\R_+\rightarrow \R$ will be called {\bf asymptotically\ uniform},
(a.u.) for short, if for any $\varepsilon>0$ there exist $T\geq 0$ and $\delta>0$
such that for all $s,t\geq T$ with $|t-s|\leq  \delta$, one has $|f(t)-f(s)|<\varepsilon$.
\end{definition}





\begin{lemma}\label{convau}
Any function $f:\R_+\rightarrow \R$ converging at infinity is (a.u.).
\end{lemma}
\begin{proof}
Let $f:\R_+\rightarrow \R$ such that $\lim_{t\rightarrow\infty}f(t)=\alpha\in \R$.
For every $\varepsilon >0$ there exists $T>0$ such that
$|f(t)-\alpha|<\frac{\varepsilon}{2}$ for $t>T$.
Hence for $s,t>T, |f(t)-f(s)|\leq|f(t)-\alpha|+|\alpha-f(s)|<\varepsilon$
independently of $|t-s|$. Thus $f$ is (a.u.).
\end{proof}

\subsubsection*{Differential case}

\begin{theorem}\label{diffcase}
Let $f:\R_+\rightarrow \R$ be differentiable with $\lim_{t\rightarrow\infty}f(t)=\alpha\in \R$.
Then $\lim_{t\rightarrow\infty}\dot f(t)=0$ if and only if
$\dot f$ is (a.u.).
\end{theorem}
\begin{proof}
According to Lemma \ref{convau}, the convergence of $\dot f$
in $\R$ entails its uniform asymptoticity.
Conversely, if $\dot f$ is (a.u.), for every $\varepsilon>0$
there exist $T\geq 0$
and $\delta>0$ such that $\forall\ s,t \geq T, |t-s|\leq \delta$\bcolor{,} we have
$|\dot f(t)-\dot f(s)|<\frac{\varepsilon}{2}$. Let $t_n=n\delta$ for $n\in\N$.
By the mean value theorem there exists
$\xi_n\in\ (t_n,t_{n+1})$ such that
$\frac{f(t_{n+1})-f(t_n)}{\delta}=\dot f(\xi_n)$. Since the left-hand side 
converges to $0$ as $n\rightarrow\infty$,
there exists $N\geq 0$ such that $|\dot f(\xi_n)|<\frac{\varepsilon}{2}$ for all $n\geq N$.
Consequently, for every $t> \max\{T,t_N\}\ \exists\  n_*\geq N$ such that
$t \in [t_{n_*},t_{n_*+1}]$; therefore, according to the triangle inequality
$$|\dot f(t)|\leq |\dot f(t)-\dot f(\xi_{n_*})|+|\dot f(\xi_{n_*})|<\varepsilon,$$
and thus $\lim_{t\rightarrow\infty}\dot f(t)=0$.
\end{proof}

We can even say more when the derivative is continuous.

\begin{theorem}\label{uniformcontinuity}
Let $f:\R_+\rightarrow \R$ be differentiable with a continuous derivative and
$\lim_{t\rightarrow\infty}f(t)=\alpha\in \R$. Then $\lim_{t\rightarrow\infty}\dot f(t)=0$ if and
only if $\dot f$ is uniformly continuous.
\end{theorem}

\begin{proof}
This follows from Theorem  \ref{diffcase} and Property 5 in Section 4.
\end{proof}

\subsubsection*{Integral case}

\begin{theorem}\label{intcasebasic}
Let $g:\R_+\rightarrow \R$ be Riemann
integrable over every interval $[0,t]$
and let $\lim_{t\rightarrow\infty}\int_0^tg(s)ds=\alpha\in \R$.
Then $\lim_{t\rightarrow\infty}g(t)=0$ if and only if $g$ is (a.u.).
\end{theorem}
\begin{proof}
Lemma \ref{convau} again entails the uniform asymptoticity of $g$ if
the latter converges. For the converse, let $\varepsilon>0$ be given.
Then there exist $T\geq 0$
and $\delta>0$ such that $\forall\ s,t \geq T, |t-s|\leq \delta$
we have
$| g(t)-g(s)|<\frac{\varepsilon}{2}$.
Since $\lim_{t\rightarrow\infty}\int_t^{t+\delta}g(s)ds = 0$, there exists
$T'\geq 0$ such that $\forall\ t>T', |\int_t^{t+\delta}g(s)ds|<\frac{\varepsilon}{2} \delta$.
According to the triangle inequality, for every $t>\max\{T,T'\}$,
$$|g(t)|\delta=\left|\int_t^{t+\delta}g(t)ds\right|\leq \int_t^{t+\delta}|g(t)-g(s)|ds
+ \left|\int_t^{t+\delta}g(s)ds\right|<2\frac{\varepsilon}{2} \delta.$$
Consequently, $|g(t)|<\varepsilon$ for $t>\max\{T,T'\}$ and thus
$\lim_{t\rightarrow\infty}g(t)=0$.
\end{proof}


As in the differential case we have the following theorem, when $g$ is continuous.

\begin{theorem}
Let $g:\R_+\rightarrow \R$ be continuous and Riemann
integrable over every interval $[0,t]$.
If $\lim_{t\rightarrow\infty}\int_0^tg(s)ds=\alpha\in \R$,
then $\lim_{t\rightarrow\infty}g(t)=0$ if and only if $g$ is uniformly continuous.
\end{theorem}

Note that the sufficient condition is the already mentioned Barb\u alat's lemma
 \cite{Far-Weg}.

\section{Some properties and examples of\\
asymptotically uniform functions}

\begin{enumerate}
\item
{\it The (a.u.) condition is purely asymptotic}, in the sense that
modifications of the function on a bounded set will not
affect the validity of the property.
\item As a direct consequence of the triangle inequality,
{\it the set of (a.u.) functions is a vector space}.
\item 
{\it If a sequence $(f_n)_{n\geq 1}$ of (a.u.) functions converges uniformly to $f$,
then $f$ is also (a.u.).} Indeed,
 for $\varepsilon>0$,
let $N$ be such that for all $n\geq N$
and $t\in \R_+$, $|f_n(t)-f(t)|<\frac{\varepsilon}{3}$. Since $f_N$ is (a.u.), by choosing $T$ sufficiently large, one can find $\delta>0$ such that for all $s,t\geq T,
|t-s|<\delta$, $|f_N(t)-f_N(s)|<\frac{\varepsilon}{3}$ and
$|f(t)-f(s)|\leq |f(t)-f_N(t)|+|f_N(t)-f_N(s)|+|f_N(s)-f(s)|<\varepsilon$.
 
\item {\it Every uniformly continuous function $f:\R_+\rightarrow \R$ is (a.u.)}: the definition of asymptotic uniformity is clearly satisfied for such a function, with $T=0$.
Consequently, since every Lipschitz continuous function
is uniformly continuous, it is (a.u.).
Moreover,
{\it any differentiable function $f:\R_+\rightarrow \R$ with a bounded derivative is (a.u.)},
since the mean value theorem entails its uniform continuity.

\item {\it A continuous function $f:\R_+\rightarrow \R$ is (a.u.) if and only if it is
uni\-form\-ly
continuous}. Indeed,
for $\varepsilon>0$, let $T$ and $\delta$ be numbers provided by
the definition of asymptotic uniformity.
Since $f$ is uniformly continuous over $[0,T+\delta]$, there exists $\delta^*$
such that for all $s,t\leq T+\delta$ with $|t-s|\leq  \delta^*$
we have $|f(t)-f(s)|<\varepsilon$.
The number
$\min\{\delta,\delta^*\}$ yields the uniform continuity over $\R_+$.
The converse is given by item 4. This result suggests that asymptotic uniformity is
mainly
interesting for asymptotically discontinuous functions. For instance,
the everywhere discontinuous function $e^{-t}I_{\Q \cap \R_+}(t)$,
which obviously converges to $0$ at infinity,  is (a.u.),
according to Lemma \ref{convau}.

\item Let $f$ and $g$ be two real-valued functions defined over $\R_+$.
{\it If $g$ is (a.u.) and $\lim_{t\rightarrow \infty}(f(t)-g(t))=0$, then $f$ is (a.u.) as well}:
indeed,
for any $\varepsilon>0$ there exist $T_1\geq 0, \delta>0$ and $T_2\geq 0$
such that $|g(t)-g(s)|<\frac{\varepsilon}{3}$ for all $s,t\geq T_1,|t-s|\leq \delta$
and $|f(t)-g(t)|<\frac{\varepsilon}{3}$ for all $t\geq T_2$; consequently, for all
$s,t\geq T=\max\{T_1,T_2\}$
with $|t-s|\leq \delta$, we have
$$|f(t)-f(s)|\leq |f(t)-g(t)|+|g(t)-g(s)|+|g(s)-f(s)|<\varepsilon.$$

\item Since a constant function is uniformly continuous, it is (a.u.).
Property 6 with $g\equiv 0$ thus confirms Lemma \ref{convau}.

\item
{\it Examples}:
\begin{enumerate}
\item The function $f(t)=t$ is unbounded but uniformly continuous and thus (a.u.). 
\item The function $f(t)=\sin(t)$, which obviously does not converge at infinity,
is uniformly continuous and thus (a.u.).
\item The bounded function $f(t)=\sin(t^2)$ is continuous
but not uniformly continuous. According to property 5,
it is not (a.u.).

\end{enumerate}

\item  
{\it An (a.u.) function is not necessarily measurable.}
Consider any boun\-ded non-measurable one,
for example the indicatrix $g$ of a non-mea\-su\-rable subset of $\R_+$.
Then $f(t)=e^{-t}g(t)$ converges to $0$ at infinity and thus is (a.u.) according to
Lemma \ref{convau}, but clearly not measurable.

\item {\it 
We now give a differentiable function which converges at infinity,
whose derivative is discontinuous on an unbounded sequence
of points, but is nevertheless
(a.u.) as it converges to $0$ at infinity (Lemma 2).}
%
Consider the function
$$
f(t) :=
\begin{cases}
1,&0\le t < 1,\\
0,&1\le t < 2,
\end{cases}
$$
periodically extended over $\R_+$.
$f$ is discontinuous at the integers, and bounded.
Clearly,
it admits a primitive, say $g$, which satisfies $0\le \dot g(t) = f(t)\le 1$.
Choosing $g(0) = 0$, we get $0\le g(t)\le t$. Then the
primitive of $e^{-t}\dot g(t)$ is $h(t) = e^{-t}g(t) - \int_0^t (-e^{-s})g(s)ds
= e^{-t}g(t) + \int_0^t e^{-s}g(s)ds$,
and the last integrand is continuous. 
Since $h$ is increasing and bounded, it
converges as $t\rightarrow \infty$. $h$ is our function: its derivative 
converges to $0$ and is discontinuous at $t=k$, $k\in \N$.

\item
Contrarily to its derivative, the unbounded function
$f(t)=t^2$ is not uniformly continuous.
Our investigations led us to the following result,
which is certainly known, but that we
did not find in the literature:
{\it if $f:\R_+\rightarrow \R$ is bounded and differentiable with a uniformly continuous
de\-ri\-va\-tive $\dot f$, then the latter is bounded and $f$ is
uniformly continuous.}
\\
Indeed, let $|f(t)|<C<\infty$ for all $t\in\R_+$. As $\dot f$ is uniformly continuous,
choosing $\varepsilon=1$, there exists $\delta>0$ such that for all $0<h<\delta$
and $t\geq 0$, $|\dot f(t+h)-\dot f(t)|<1$. Furthermore, for any fixed $h_0$ with
$0<h_0<\delta$ and $t\geq 0$, according to the mean value theorem, there exists
$\xi(t) \in (t,t+h_0)$ such that $f(t+h_0)-f(t)=\dot f(\xi(t))h_0$. We thus get
$|\dot f(\xi(t))|=|\frac{f(t+h_0)-f(t)}{h_0}|\leq\frac{2C}{h_0}$ and for every $t\in R_+$
$$|\dot f(t)|\leq |\dot f(t)-\dot f(\xi(t))|+|\dot f(\xi(t))|< 1+\frac{2C}{h_0}.\qquad\qquad\qquad\QEDB$$

\item It follows from the definition that
{\it if $f:\R_+\rightarrow \R$  is (a.u.) and if $g:\R\rightarrow \R$ is uniformly
continuous, then the composite function $g\circ f$ is (a.u.) as well}.

\item {\it The (a.u.) property is not preserved by almost everywhere equality.}
For instance, $f(t)=e^{-t}$ is (a.u.) but clearly not $g(t)=e^{-t}+I_{\N}(t)$,
despite the fact that the two functions are equal almost everywhere.

\end{enumerate}

\section{About converging functions}

\begin{definition}
We shall say that $f: \R_+\to \R$ satisfies condition $(C)$, if
$$\textrm{for a.e.\ (i.e., \ for \ almost \ every}) \ t \in \R_+ ,\ 
\lim_{n\rightarrow \infty}f(nt)=0\leqno(C)$$
and condition $(C^*)$ if
$$\textrm{for all}\ t \in \R_+ ,\ 
\lim_{n\rightarrow \infty}f(nt)=0.\leqno(C^*)$$
\end{definition}

Condition (C) is for instance satisfied for Lebesgue integrable functions,
as shown in \cite{Les}.
The following result holds for continuous functions.

\begin{lemma}[Croft's lemma {\cite[p.\ 17]{And-Mor-Tet}}]\label{Croft}
If $f: \R_+ \rightarrow \R$ is continuous and $(C^*)$ holds, then
$ \displaystyle \lim_{t\rightarrow \infty} f(t)=0$.
\end{lemma}

On the other hand,
$(C)$ is sufficient to provide convergence for uniformly continuous functions.
We prove the following equivalence.

\begin{lemma}\label{Cunifcont}
Let $f: \R_+ \rightarrow \R$ be continuous and $(C)$ hold.
Then the following two sta\-te\-ments are equivalent:
\begin{enumerate}[(a)]
\item
$\displaystyle \lim_{t\rightarrow \infty}f(t)=0$;
\item
$f$ is uniformly continuous.
\end{enumerate}
\end{lemma}

\begin{proof}
%
%
%
(a) $\Rightarrow$ (b): if a continuous function converges at infinity,
it is uniformly continuous \cite[p.\ 233]{And-Mor-Tet}.

\noindent
(b) $\Rightarrow$ (a):
let $E$ be the set of all $t \in \R_+$ for which $f(nt)$ does not converge to $0$ as
$n\uparrow \infty$. Suppose that $E\neq \emptyset$ and let $t\in E$.
By definition of the convergence to 0,
there exist $\varepsilon_0>0$ and $(n_k)_{k\geq 1} \uparrow \infty$
such that
$$\forall\ \ k\geq 1,\qquad |f(n_kt)|>\varepsilon_0.$$

According to the uniform continuity of $f$,
$$\forall\ \ \varepsilon>0, \exists \ \delta (\varepsilon) \text{ such that } \forall\ \ t,s\geq 0,
|t-s|\leq \delta(\varepsilon)\Rightarrow |f(t)-f(s)|<\varepsilon.$$

Let $\varepsilon'=\frac{\varepsilon_0}{2}$.
If $\delta(\varepsilon')\in E^c$, we choose $\delta' = \delta(\varepsilon')$.
Since the Lebesgue measure of $E$ is zero, $E^c$ is dense in $\R_+$.
If $\delta(\varepsilon')\in E$, it is thus always pos\-sible to choose
$\delta' \in E^c$ such that $0<\delta'<\delta(\varepsilon')$. 
It follows that
the sequence $(f(n\delta'))_{n\geq 1}$ converges to $0$ as $n\uparrow \infty$.
Therefore, there exists $N\geq 1$ such that

$$\forall\ \ n\geq N,\qquad |f(n\delta')|<\varepsilon'=\frac{\varepsilon_0}{2}.$$
We partition $\R_+$ according to
$$\R_+=\bigcup_{n\geq 0}[n\delta', (n+1)\delta')
$$
and choose $n_{k_0}>N$ and $n'>N$ such that $n_{k_0}t$ belongs to the interval 
$[n'\delta', (n'+1)\delta')$. Then
$$\varepsilon_0<|f(n_{k_0}t)|\leq |f(n_{k_0}t)-f(n'\delta')|+|f(n'\delta')|<2\frac{\varepsilon_0}{2}=\varepsilon_0.$$
This contradiction entails $E=\emptyset$ and thus $(C^*)$ holds. $f$ being continuous,
Lemma \ref{Croft} (Croft's lemma) completes the proof.
\end{proof}


One obtains another equivalence when the continuity assumption is lifted.

\begin{theorem}\label{Cuplusr}
Let $f: \R_+ \rightarrow \R$ be such that $(C)$ holds. 
The following two sta\-te\-ments
are equivalent:

\begin{enumerate}[(a)]

\item $\displaystyle \lim_{t\rightarrow \infty}f(t)=0$;

\item $f=u+r$ with $u$ uniformly continuous and $\lim_{t\rightarrow \infty}r(t)=0$.

\end{enumerate}
\end{theorem}

\begin{proof}
(a) $\Rightarrow$ (b) is obvious, since we can choose $u(t)\equiv 0$ and $f=r$.
\\
\noindent
(b) $\Rightarrow$ (a):  since (C) holds, for a.e.\ $t>0$
$$\lim_{n\rightarrow \infty}(u(nt)+r(nt))=0$$
and because $r$ converges to $0$, for a.e.\ $t>0$,
 $$\lim_{n\rightarrow \infty}u(nt)=0.$$
Thus $(C)$ is also valid for $u$. Since the latter is uniformly continuous, Lemma
\ref{Cunifcont}
provides its convergence to zero at infinity, and thus that of $f$.
\end{proof}

\begin{definition}
We shall say that $f: \R_+ \rightarrow \R$ is a $(u,r)$--function if
$f=u+r$, with $u$ uniformly continuous and $r$ converging to $0$ at infinity.
\end{definition}

{\em Remark\ } 
This decomposition is
obviously not unique: for instance, one may add to $u$ any uniformly continuous function
vanishing at infinity and subtract it from $r$.

\section{Back to asymptotically uniform functions}

Let $\N^*$ denote the positive integers.
In our context, a function $f:\R_+\rightarrow \R$ is {\it piecewise differentiable}
if there exists a strictly increasing sequence $(t_n)_{\in \N^*}$ going to $\infty$ as $n\rightarrow \infty$ such that $f$ is required to be everywhere continuous and
differentiable only over the open intervals $(t_n,t_{n+1}), n\in \N^*$.
The following consequence of the mean value theorem is likely
known. We nevertheless provide a proof, which
we did not find in the literature.

\begin{lemma}\label{piecdifflipschitz}
Let $f:\R_+\rightarrow \R$ be piecewise differentiable with associated sequence
$(t_n)_{n\in \N^*}$. If $|\dot f(t)|\leq M< \infty$ over $\bigcup_{n\geq 1}(t_n,t_{n+1})$, 
then $f$ is Lipschitz continuous.
\end{lemma}
\begin{proof}
Let $0\leq s< t$. If the interval $[s,t]$ does not contain elements of $(t_n)_{n\in \N}$,
the mean value theorem entails $|f(t)-f(s)|\leq M(t-s)$.
Otherwise there exist $n_p\leq n_q$ in $\N^*$ such that
$0\leq s \leq t_{n_p}< t_{n_p+1}< ...<t_{n_q} \leq t,$
where $t_{n_p}$ is the smallest element of the sequence $(t_n)_{n\in \N^*}$
to the right of $s$ and
similarly for $t_{n_q}$ to the left of $t$ ($n_p= n_q$ means that $[s,t]$ contains only
one element $t_{n_p}$ with $0\leq s\leq  t_{n_p}\leq t)$ .
The triangle inequality implies
$$|f(t)-f(s)|\leq |f(t_{n_p})-f(s)|+|f(t_{n_{p}+1})-f(t_{n_p})|+...+|f(t)-f(t_{n_q})|.$$
Since $f$ is continuous over each interval $[t_n,t_{n+1}]$ and differentiable
over $(t_n,t_{n+1})$, the mean value theorem applies and 
$$|f(t)-f(s)|\leq M(t_{n_p}-s)+M(t_{n_{p}+1}-t_{n_p})+...+M(t-t_{n_q})= M(t-s).$$
Hence $f$ is Lipschitz continuous with constant $M$.
\end{proof}


For $f:\R_+\rightarrow \R$, $T\geq 0$ and $\varepsilon >0$,
we define
$V^{\varepsilon}_T(f):=\bigcup_{t\geq T}(f(t)-\varepsilon,f(t)+\varepsilon)$.

\begin{theorem}\label{aulipschitz}
For $f:\R_+\rightarrow \R$,
the following two statements are equivalent:
\begin{enumerate}[(a)]
\item
$f$ is (a.u.);
\item
$\forall\ \varepsilon>0, \exists \ T\geq 0$ such that $V^{\varepsilon}_T(f)$ contains
the graph of a function defined over $[T,\infty)\,$ and Lipschitz continuous there.
\end{enumerate}
\end{theorem}
\begin{proof}
(a) $\Rightarrow$ (b):
for $\varepsilon>0$, let
$T\geq 0$ and $\delta>0$
be numbers guaranteed by the definition of the uniform asymptoticity of $f$
for $\frac{\varepsilon}{2}$.
Consider the equispaced nodes (points at which a function is evaluated)
$t_n=T+n\frac{\delta}{2}$ for all $n\in\N$, let $g(t_n)=f(t_n)$ and
define $g$ over each
interval $[t_n,t_{n+1}]$ as the segment connecting $\big(t_n,g(t_n)\big)$
and $\big(t_{n+1},g(t_{n+1})\big)$
($g$ is clearly continuous over $[T,\infty)$).
Since $t_{n+1}-t_n=\frac{\delta}{2}$, the uniform asymptoticity of $f$ entails
$|f(t)-f(t_n)|<\frac{\varepsilon}{2}$ for all $t \in[t_n,t_{n+1}]$. Moreover, $g$ being
affine 
over $[t_n,t_{n+1}]$,
$$|f(t_n)-g(t)| = |g(t_n)-g(t)|
= |a_n||t_n-t| \le |a_n||t_n - t_{n+1}| = |g(t_n)-g(t_{n+1})|,$$
where $a_n$ is the slope of $g$
over $[t_n,t_{n+1}]$ (a point moving on a line always drifts away from
its initial vertical position).
Thus $|f(t_n)-g(t)| \leq |g(t_n)-g(t_{n+1})| = |f(t_n)-f(t_{n+1})|<\frac\varepsilon 2$
for all $t$ in $[t_n,t_{n+1}]$.
Accordingly, for all  $n\in \N$, we get
$$|f(t)-g(t)|\leq |f(t)-f(t_n)|+|f(t_n)-g(t)|<2\frac{\varepsilon}{2}=\varepsilon$$
for all $t$ in $[t_n,t_{n+1}]$. Thus $|f(t)-g(t)|< \varepsilon$ for all $t\in [T,\infty)$
and the graph of $g$ belongs to $V^{\varepsilon}_T(f)$.
Furthermore, the slope
of $g$ over $[t_n,t_{n+1}]$ satisfies
$$|a_n|=\bigg|\frac{f(t_{n+1})-f(t_n)}{\frac{\delta}{2}}\bigg|<\frac{\frac{\varepsilon}{2}}
{\frac{\delta}{2}}=\frac{\varepsilon}{\delta},$$
a bound independent of $n$.
$g$ being continuous over $[T,\infty)$ and piecewise differentiable with bounded derivative, it is Lipschitz continuous according to Lemma \ref{piecdifflipschitz}.
\\
(b) $\Rightarrow$ (a):
let $T$ correspond to $\frac{\varepsilon}{3}>0$
in point (b) and let $g$ be a 
Lipschitz function whose graph is contained in $V^{\frac{\varepsilon}{3}}_T(f)$.
As $g$ is obviously uniformly continuous,
$\exists \ \delta>0$ such that $\forall\ s,t\geq T, |t-s|<\delta$,
we have $|g(t)-g(s)|<\frac{\varepsilon}{3}$ and thus
$$|f(t)-f(s)|\leq |f(t)-g(t)|+|g(t)-g(s)|+|g(s)-f(s)|<3\frac{\varepsilon}{3}=\varepsilon.$$
Consequently, $f$ is (a.u.).
\end{proof}

\begin{corollary}
For a function $f:\R_+\rightarrow \R$, the following three statements are equivalent:

\begin{enumerate}[(a)]
\item
$f$ is (a.u.);
\item
There exists a uniformly continuous function $u:\R_+\rightarrow \R$ such that
$\lim_{t\rightarrow \infty}(f(t)-u(t))=0$\bcolor{;}
\item
$f$ is a $(u,r)$--function.
\end{enumerate}
\end{corollary}

\begin{proof}
(a) $\Rightarrow$ (b):
we shall construct $u$ as a piecewise affine function\bcolor{,} as in the first part
of the proof of Theorem \ref{aulipschitz}, but with a graph contained in strips of
decreasing thickness.
For $\varepsilon>0$ \ and every $k\geq 1$
we define $\varepsilon_k=\frac{\varepsilon}{k}$ and denote by $u_k$ the function
corresponding to $\frac{\varepsilon_k}{2}$ and $t\geq T_k$ in the first part
of the proof of Theorem \ref{aulipschitz}. 
As we want $T_k$ to be a node,
i.e., $u_{k-1}(T_k)=f(T_k)=u_k(T_k)$,
we choose the corresponding $\delta_k$ as $\delta_{k-1}$ divided by some integer
power of $2$. Furthermore, we choose $T_k$ sufficiently large to warrant
that 
$T_k>T_{k-1}+1$.
Now define the function $u$
to be $u_k$ over $[T_k,T_{k+1}]$ and $u_1(T_1)$ over $[0,T_1]$.
Clearly $u$ is continuous and $|f(t)-u(t)|<\varepsilon_k$ for $t\geq T_k$;
since $\varepsilon_k\downarrow 0$ as 
$k \rightarrow \infty$, we deduce that  $\lim_{t\rightarrow \infty}(f(t)-u(t))=0$.
According to property 6, $u$ is (a.u.). Since $u$ is continuous, property 5 entails its uniform continuity and we get (b).\\
(b) $\Rightarrow$ (a):
by property 4, $u$ is (a.u.), and by property 6,
$f$ is (a.u.).\\
Since (b) $\Leftrightarrow$ (c) is obvious, we have proved the equivalence of the three
statements.
\end{proof}

\vfill\eject
\section{Extension to higher-order derivatives}

\begin{theorem}
Let $f:\R_+\rightarrow \R$ be $n$-times differentiable
with
$\lim_{t\rightarrow\infty}f(t)=\alpha \in \R$. Then the following two assertions are equivalent:
\begin{enumerate}[(a)]
\item
$\forall\ \ 1\leq k\leq n, \ \lim_{t\rightarrow\infty}f^{(k)}(t)=0$;
\item
$f^{(n)}$ is (a.u.).
\end{enumerate}
\end{theorem}

Note that the case $n=1$ is Theorem \ref{diffcase}.

\begin{proof}
(a) $\Rightarrow$ (b):
Since $\lim_{t\rightarrow\infty}f^{(n)}(t)=0$, 
$f^{(n)}$ is (a.u.) by Lemma 2.
\\
(b) $\Rightarrow$ (a):
We shall say that 
a function $F$ defined over $\R_+$ ultimately satisfies a property $\mathbb P$
if there exists $T\geq 0$ such that $F$ has the property $\mathbb P$ for all $t\geq T$.

To prove (a), it is sufficient
to show that $f^{(n)}$ is ultimately
bounded.
Indeed, if this is true, the already mentioned Theorem 6 in \cite[p.\ 141]{Cop}
provides the convergence to zero of $f^{(k)}(t)$,
for $1\leq k\leq n-1$, when $t\rightarrow\infty$.
Since the derivative of $f^{(n-1)}$ is (a.u.) by hypothesis,
the convergence of
$f^{(n)}$ to $0$ is given by that of $f^{(n-1)}$
and Theorem \ref{diffcase}.

First we note that, under the hypothesis,
none of the derivatives is ultimately $0$.
Indeed, if one of them is 0,
let $i, 1\leq i\leq n,$ be the smallest index such that
$f^{(i)}$ is identically $0$ over some interval $[T, \infty)$.
If $i=1$, then $f$ is constant over $[T, \infty)$ and
(a) is clearly true since all subsequent derivatives are also $0$
over $(T, \infty)$.
On the other hand, under the hypothesis,
the case $i\geq 2$ leads to a con\-tra\-dic\-tion.
Indeed, the first derivative would have to be
different from $0$ asymptotically; moreover,
$f$ would be a polynomial of degree $\geq 1$.
In that case $f$ could not converge
in $\R$ as $t\rightarrow \infty$.
We can thus assume that none of the derivatives is ultimately $0$.

Let $f:\R_+\rightarrow \R, n\in \N^*, t \geq 0, h>0$.
According to Taylor's theorem over closed intervals \cite[p.\ 168]{Fle},
if $f^{(n-1)}$ is continuous over $[t,t+h]$ and $f^{(n)}$
defined over $(t,t+h)$, then there exists $t<\xi_{t,h}<t+h$ such that
$$f(t+h)=f(t)+\sum_{k=1}^{n-1}f^{(k)}(t)\frac{h^k}{k!}+f^{(n)}(\xi_{t,h})\frac{h^n}{n!}.$$
$f$ being $n$ times differentiable over $\R_+$, it satisfies these
conditions and thus
\begin{align*}
\sum_{k=1}^{n-1}f^{(k)}(t)\frac{h^k}{k!}&=f(t+h)-f(t)-f^{(n)}(\xi_{t,h})\frac{h^n}{n!}\\
&=f(t+h)-f(t)-(f^{(n)}(\xi_{t,h})-f^{(n)}(t))\frac{h^n}{n!}-f^{(n)}(t)\frac{h^n}{n!}.
\end{align*}
Extending the index range to $n$, we get
$$R_h(t):=\sum_{k=1}^{n}f^{(k)}(t)\frac{h^k}{k!}=f(t+h)-f(t)-(f^{(n)}
(\xi_{t,h})-f^{(n)}(t))\frac{h^n}{n!}.$$
Since $f$ converges to $\alpha \in \R$ when $t\rightarrow\infty$ and $f^{(n)}$ is (a.u.),
for every $M>0$, there exists $T_{M}\geq 0$ and $0<\delta_{M}$
(we can always choose $\delta_{M}<1$), so
that $\forall\ t\geq T_{M}$ and
$0<h\leq \delta_{M}<1$
$$|f(t+h)-f(t)|<\frac{M}{2} \ \textrm{and} \ |f^{(n)}(\xi_{t,h})-f^{(n)}(t)|\leq \frac{M}{2}.$$
Thus
\begin{align*}
\left|(f^{(n)}(\xi_{t,h})-f^{(n)}(t))\frac{h^n}{n!}\right|&=\left|f^{(n)}(\xi_{t,h})-f^{(n)}(t)
\right|\frac{h^n}{n!}\\
&\leq \left|f^{(n)}(\xi_{t,h})-f^{(n)}(t)\right|\leq \frac{M}{2}.
\end{align*}
Hence, for every $h$ with $0<h\leq \delta_{M}<1$, $|R_h(t)|\leq M$
$\forall\ \ t \in[T_{M},\infty)$.
Consequently 
$R_h$ is ultimately bounded, uniformly
in $h$ for $0<h\leq \delta_{M}<1$.
Clearly, if two functions are bounded over $[T_{M},\infty)$,
then the same property holds for any multiple or linear combination of them.
In the spirit of Richardson extrapolation \cite[p.\ 292]{Hil},
we shall use the two operations to successively eliminate
all the terms of $R_h$ except the last one,
which will be a bounded multiple of $f^{(n)}$; $f^{(n)}$ will thus bounded
over $[T_{M},\infty)$ as well.
We introduce the functions $f_k(t):=\frac{f^{(k)}(t)}{k!}$
and recall that none of them is ultimately $0$. From the definition we get
\begin{align*}
R_{h}&=f_1h+f_2h^2+...+f_nh^n,\\
R_{\frac{h}{2}}
&=f_1\frac{h}{2}+f_2\frac{h^2}{2^2}+...+f_n\frac{h^n}{2^n}.
\end{align*}
Then one consider\bcolor{s} the function
$$R^{(1)}_{h}:=R_{h}-2R_{\frac{h}{2}}=0+f_2(1-\frac{1}{2})h^2
+f_3(1-\frac{1}{2^2})h^3+...+f_n(1-\frac{1}{2^{n-1}})h^n,$$
which,
by the triangle inequality, satisfies $|R^{(1)}_{h}|\leq (1+2)M$ over
$[T_{M},\infty)$. The next step yields
\begin{align*}
R^{(2)}_{h}:=&\ R^{(1)}_{h}-2^2R^{(1)}_{\frac{h}{2}}\\
=&\ 0
+f_3(1-\frac{1}{2^2})(1-\frac{1}{2})h^3+...
+f_n\big(1-\frac{1}{2^{n-1}}\big)\big(1-\frac{1}{2^{n-2}}\big)h^n;
\end{align*}
again, by the triangle inequality, $|R^{(2)}_{h}|\leq (1+2^2)(1+2)M$.
After $n-1$ steps, we get a single term,
which is a multiple of $f^{(n)}$, namely:
$$R^{(n-1)}_h=R^{(n-2)}_h-2^{n-1}R^{(n-2)}_{\frac{h}{2}}
=\big(1-\frac{1}{2^{n-1}}\big)\big(1-\frac{1}{2^{n-2}}\big)
\cdot \cdot \cdot \frac{1}{2}\frac{h^n}{n!}f^{(n)}.$$
Since $|R^{(n-1)}_h|\leq (1+2^{n-1})(1+2^{n-2})\cdot \cdot \cdot(1+2)M$,
$f^{(n)}$ is bounded over $[T_{M},\infty)$.
\end{proof}

\vskip1em
We can also adapt Hadamard's lemma
(Theorem \ref{th:Hadamard}) to the case when the derivatives are
uniformly continuous. This also generalizes 
Theorem \ref{uniformcontinuity}.

\begin{theorem}
Let $f:\R_+\rightarrow\R$ be $n$-times differentiable with $f^{(n)}$ continuous
and assume that $\lim_{t\rightarrow\infty}f(t) = \alpha \in \R$.
The following two statements are equivalent:
\begin{enumerate}[(a)]
\item
All the derivatives $f^{(1)},\ldots,f^{(n)}$ converge to $0$ at infinity;
\item
All the derivatives $f^{(1)},\ldots,f^{(n)}$ are uniformly 
continuous.
\end{enumerate}
\end{theorem}
We stress that the result holds for one more derivative than
in Hadamard's lemma.

%
%
 \bibliographystyle{}
 \bibliography{alpha}
 \def\backskip{\hskip-0.35pt}
\def\book #1#2#3#4#5{\par\vskip0.2em\noindent\hangindent3.5em#1\backskip #2:
{\it #3}. #4 (#5)}
\def\art #1#2#3#4#5#6#7{\par\vskip0.2em\noindent\hangindent3.5em
#1\backskip #2: #3. {\it #4} {\bf #5}, #7 (#6)}
\def\contr #1#2#3#4#5#6#7#8{\par\vskip0.2em\noindent\hangindent3.5em
#1\backskip #2: #3. In: #4, {\it #5}, pp.\ #7. #6 (#8)}
\def\contrsaut #1#2#3#4#5#6{\par\vskip0.2em\noindent\hangindent3.5em
#1 #2, #3, in: {\it #4} (#5) #6}
\def\thiss #1#2#3#4{\par\vskip0.2em\noindent\hangindent3.5em
#1 #2, #3, {\it #4}}
\def\preprint #1#2#3#4{\par\vskip0.2em\noindent\hangindent3.5em#1\backskip #2:
{\it #3}, #4}
\def\toappear #1#2#3#4{\par\vskip0.2em\noindent\hangindent3.5em#1 #2,
#3, {to appear in \it #4}}
\def\vk #1#2#3{\par\vskip0.8em\noindent{\bf\parindent2em#1.\ #2\hfill#3}}
\def\va #1#2#3{\par\vskip0.4em\noindent{\rm\hskip1em\parindent
4em#1.\ #2\hfill#3}}
\def\vua #1#2#3{\par\vskip0.2em\noindent\hskip2em\parindent
6em#1.\ #2\hfill#3}
\def\vuua #1#2#3{\par\noindent\hskip3em\parindent
8em#1.\ #2\hfill#3}

\end{document}